**Distribution Agreement**

Signature:

______________________________ _______________

Dalton Bidleman Date

Riemann—Roch for Toric Rank Functions

By

Dalton Bidleman
Master's of Science

Mathematics

David Zureick Brown, Ph.D.
Advisor

Ron Gould, Ph.D.
Committee Member

Shawn Ramirez, Ph.D.
Committee Member

Accepted:

Lisa A. Tedesco, Ph.D.
Dean of the Graduate School

Date

Riemann—Roch for Toric Rank Functions

By

Dalton Bidleman

Master's of Science
B.S., Emory University 2015

Advisor: David Zureick-Brown, Ph.D.

An abstract of
A thesis submitted to the Faculty of the James T.
Laney School of Graduate Studies of Emory University
in partial fulfillment of the requirements for the degree of
Master of Science
in Mathematics
2015

**Abstract**

Riemann--Roch for Toric Rank Functions
By Dalton Bidleman

In this thesis we study toric rank functions for chip firing games and prove special cases of a conjectural Riemann-Roch. The original motivation for an investigation into this area of study came for the adaptation (due to Matt Baker) of Riemann-Roch into a graph theoretic analogue through the use of chip-firing games. Here, we collect known results and present new observations that indicate Riemann--Roch holds for trees and polygons. We also prove an asymptotic case of Riemann--Roch (i.e.~Riemann--Roch for divisors of large degree). Finally, we also provide magma code and computational evidence that Riemann--Roch holds for the toric rank function.

Riemann—Roch for Toric Rank Functions

By

Dalton Bidleman

Master's of Science
B.S., Emory University 2015

Advisor: David Zureick-Brown, Ph.D.

A thesis submitted to the Faculty of the James T.
Laney School of Graduate Studies of Emory University
in partial fulfillment of the requirements for the degree of
Master of Science
in Mathematics
2015

# Table of Contents

# List of Figures

# Chapter 1

# Introduction

The history of the Riemann–Roch Theorem traces back to 1857 when Riemann proved the initial result for complex analysis and algebraic geometry, which he called the Riemann Inequality. His student Gustav Roch finalized it in 1865. This theorem, in a general sense, says that for a surface, functions with prescribed zeroes and poles satisfy strong numerical constraints (see Chapter 2 for more precise definitions). Also, since the middle of the 20th century, there have been several papers devoted to chip-firing games [4,5,6,7,8] played on the vertices of a graph.

These games became of interest to Matt Baker and Sergei Nourine [1] when they were attempting to find a graph theoretic analogue of the Riemann–Roch Theorem. Through a rather clever combinatorial approach, they were able to show that the classical Riemann–Roch Theorem applies to graphs with slight modifications to some of the definitions. Therefore, for any graph we were now able to define the notion of the rank of a divisor on a graph in terms of other features that describe the graph. A divisor on a graph is just a labeling of the vertices, and Baker's notion of rank is a numerical measure of the "robustness" of this labeling.

This application of the theorem opened many more current doors of study that motivate the findings in this paper. Specifically, the work of Katz and Zureick-Brown [2] in 2012 utilized these ideas of rank functions, but applied them in an algebro-geometric context to study Diophantine problems, allowing them to combinatorially bound numbers of solutions to equations. Likewise,

Matt Baker was able to take his own work a step farther and prove that the Brill–Noether Conjecture has an analogue for graphs, and other authors have used this to give combinatorial proofs of the Brill–Noether Theorem [3]. This conjecture bounds the rank given certain facts about the graph hold. It is because of these things that attempts have been made to see how the Riemann–Roch Theorem acts on other surfaces and in other cases, mainly by varying variables like the genus of the graph.

This line of thinking led to the analysis of how Riemann–Roch Theorem can be applied to graph curves. Graph curves have the useful property that for every graph curve we can associate to it a dual graph (see Figure 2.2). This allows us to play a chip firing game on the dual graph. That game immediately gives us a scenario in which the Riemann–Roch Theorem is potentially applicable again. However, this new instance of chip firing is not as a straight forward because we now have to deal with the intersection points of the graph curves. That is how we formulate the idea of the toric rank that defines the whole of Chapter 3. In the first part of the chapter we define necessary terms that indicate what the individual pieces of our divisor look like, where a divisor is just some initial configuration. We also define the process of chip firing on these new structures. After that, we formalize the process by which we are able to make the intersection points agree and that is how we arrive at the definition of toric rank. It is with this new definition, that we can approach new problems that deal with the Riemann–Roch Theorem. With toric rank, we are able to prove Riemann–Roch holds for trees, for $n$-gons, and even the asymptotic case of Riemann–Roch (i.e. for divisors of large degree). To finish, in the final chapter we provide magma code that allow us to compile computational evidence toward the validity of Riemann–Roch on intermediate cases of toric rank functions. These functions can also be used more simply to look at ranks, as Matt Baker defined them [1], on graphs and their ability to meet the criterion of the graph theoretic Riemann–Roch Theorem.

# Chapter 2

# Background Information

## 2.1 Divisors on graphs

To begin, we define a graph and the associated notions that will be used to investigate toric rank functions and chip firing games.

**Definition 2.1.1.** Let $V$ be the set of vertices and $E$ the set of edges. Then, generally we can consider a graph to be an ordered pair of disjoint sets $(V, E)$ such that $E$ is a subset of the set of unordered pairs of $V$.

Unless stated, all graphs that will be dealt with in this paper will be finite; therefore $V$ and $E$ will be finite as well. Also, a graph for the sake of the arguments presented here, will always be connected and have no loops. For notation, note that $V(G)$ and $E(G)$ will simply represent the vertex and edge set on some specific graph $G$. Now, take some arbitrary edge $(x, y)$. It is said to join the vertices $x$ and $y$ and is denoted by $xy$. Therefore, it is important to note that $xy$ and $yx$ are the same edge and if $x, y \in E(G)$ then $x$ and $y$ are adjacent vertices of G and $x$ and $y$ are incident with an edge $xy$. Note that Figure 2.0 below highlights a graph with $V(G) = \{1, 2, 3, 4\}$ and $E(G) = \{12, 23, 34, 41\}$.

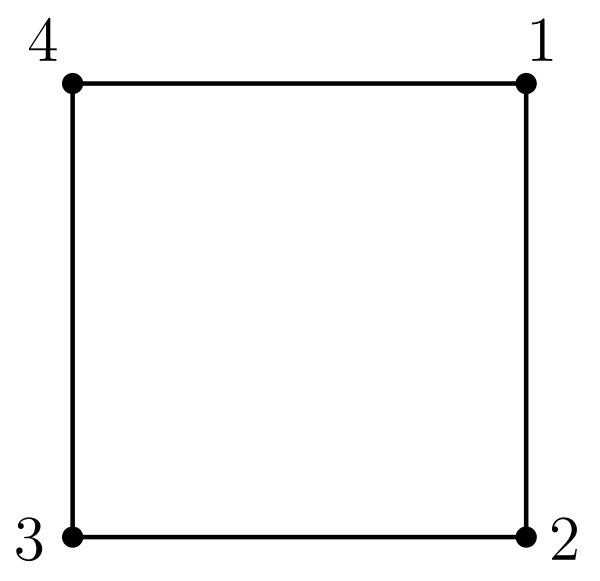

Figure 2.0 Graph with 4 vertices and 4 edges

With each graph there is this linear algebraic notion of an adjacency matrix, which acts as a bookkeeping device for a graph.

**Definition 2.1.2.** An adjacency matrix on our finite graph $G$ of $n$ vertices is an $n \times n$ matrix whose non-diagonal entries are exactly the number of edges from vertex $i$ to vertex $j$ and whose diagonal entries, are generally $0$ in the examples of this paper, but more generally represent the number of loops from a vertex $v_i$ to itself.

Therefore, recalling the graph in Figure 2.0. we see that the adjacency matrix is just

$$\begin{pmatrix} 0 & 1 & 0 & 1 \\ 1 & 0 & 1 & 0 \\ 0 & 1 & 0 & 1 \\ 1 & 0 & 1 & 0 \end{pmatrix}$$

Using the information generated by the adjacency matrix we can define the next useful notion in graph theory called a Laplacian matrix if we first make use of the degree matrix $D$. The degree matrix allows us to keep track of the number of edges that are attached to each vertex $v_i$, which is also denoted as deg$(v_i)$. The matrix itself has 0 at all non-diagonal entries and $\deg(v_i)$ at each $a_{ij}$ where $i = j$. It is from this that we also get the degree of the entire graph $G$, as it is equal to the degree of greatest value on $V(G)$. Then, the Laplacian matrix of $G$ is defined as the difference of the degree matrix $D$ and the adjacency matrix, that is, $L = D - A$. The Laplacian allows us to measure by what extent a graph differs at one vertex from the values of the ones nearby. It is clear then, that the laplacian of $G$ just takes the negative values of the adjacency matrix at all non-diagonal entries and keeps the degree matrix entries along the diagonal. Therefore, referring back to the graph in figure 2.1 we see that the Laplacian of this graph is just:

$$\begin{pmatrix} 2 & -1 & 0 & -1 \\ -1 & 2 & -1 & 0 \\ 0 & -1 & 2 & -1 \\ -1 & 0 & -1 & 2 \end{pmatrix}$$

**Definition 2.1.3.** Mathematically, the genus is the number $g$, where $g = |edges| - |vertices| + 1$. The genus is also representative topologically as the number of the holes in the graph. For clarity, both a torus and a mug with a handle have genus 1.

After this, it is only natural to look at divisors on a graph $G$.

**Definition 2.1.4.** For a graph G, the group of divisors, denoted as $Div(G)$, is the free abelian group on $V(G)$. Then, a divisor is just an integer labeling of the vertices where a divisor $D \in Div(G)$ takes the form $\sum a_v(v)$, such that the coefficients are integers that represent the value of the given vertex $v$.

**Definition 2.1.5.** Following this, a canonical divisor is defined as $K_G = \sum (deg(v) - 2)(v)$. We simply define the value at the vertex as two less the original degree. The most important fact about a canonical divisor deals with its' degree. Specifically, we know that deg$(K_G) = 2g - 2$, where $g$ is simply the genus as defined before. This is the case since the sum over all vertices $v$ of $\deg(v)$ is equal to twice the number of edges in $G$. Therefore, we get deg$(K_G)$= $2|E(G)| - 2|V(G)| = 2g - 2$ as desired.

## 2.2 Chip Firing Games

Next, we must define the notion of a chip firing game on a graph, call it G. The game is played on $V(G)$, and fixes a certain number of chips, either positive or negative, on each vertex. We have several definitions on the chip-firing game, whose graph theoretic analogue has been defined in the prior section.

**Definition 2.2.1.** The initial configuration, or a divisor $D$, in this case, has the sum of the chips in $D$ denoted as the degree of $D$, or $\deg(D)$ for simplicity.

This set up is then represented by an $n \times 1$ vector whose entries are just the chip value of that vertex on our graph $G$. We introduce the term in-debt to represent any vertex at which there are some negative number of dollars present. The goal of a chip firing game is to then manipulate the chips in such a way that every vertex is out of debt. This can be done in two ways, denoted as the two possible moves in a chip firing game. The first move is called borrowing. Borrowing consists of a vertex taking one chip from each vertex adjoined to it by an edge. For a concrete example, Figure 2.1 below highlights an instance of borrowing at the topmost vertex.

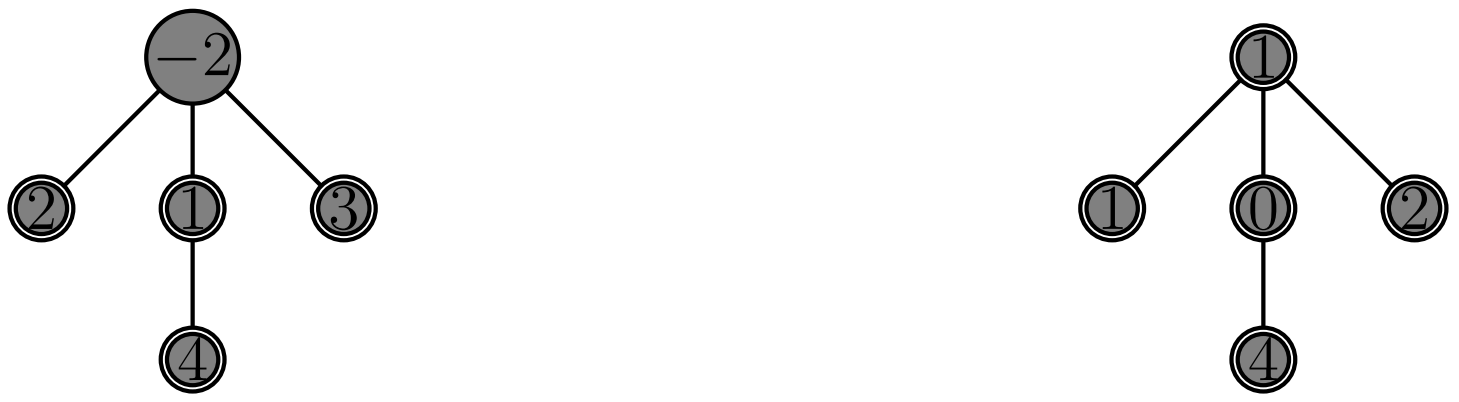

Figure 2.1: Instance of a Borrowing Move.

---

Likewise, a lending move consists of a vertex giving one chip to each vertex adjoined to it by an edge. Both of these instances are also represented as a vector, specifically the firing vector, that keeps up with how the chips are transmitted about the vertices with either a positive value representing borrowing or a negative value representing firing.

**Definition 2.2.2.** Every possibility that is achieved by some combination of lending and borrowing moves is called the divisor group $|D|$ and some arbitrary new divisor in this group is denoted as $D'$.

Furthermore, these moves can be performed multiple times and in tandem to potentially produce a graph that has every vertex out of debt. Such a configuration is called a winning position or an effective divisor and the moves it took to produce such a winning position are called a winning strategy. Likewise, a divisor in which there is no possible way to get all the vertices out of debt is called non-effective. With all the necessary definitions in place, we can use the genus as defined before to establish some facts about chip-firing games. If $\deg(D) \geq g$, then there is always a winning strategy. Also, if $\deg(D) \leq g - 1$ there is an initial configuration for which no winning strategy exists. Next, we need to define an important notion of rank for a chip-firing game.

**Definition 2.2.3.** Rank is defined as the largest numbers of chips that can be arbitrarily removed from some graph $G$ and it still have some possible winning strategy.

As a byproduct of this, a game in which the initial configuration is not winnable has $rank = -1$. An algorithmic way to calculate rank is found in a combinatorial approach to the possible way of distributing the chips. For each possible value of the rank, $r$, we remove $r$ chips arbitrarily. Therefore, for $r = 2$, we would have to check the cases when we remove 2 from some vertex and also the cases when we remove 1 from two distinct vertices. We can then write each of these possibilities as a vector, call the vector $E_i$, and check that for some $D'$ in $|D|$, $D' - E_i \geq 0$, which implies that there is a winning strategy. The instant that we find some $E_i$ whose difference with every $D' \in |D|$ is not effective, we conclude that it has rank equal to one less the sum of the entries in E. Mathematically, this can be calculated by the Riemann–Roch analog for graph theory. This says that given some divisor $D$, $r(D) - r(K - D) = \deg(D) + 1 - g$. In this formula, we have already defined all the pieces. $r(D)$ is the rank of the divisor, $r(K - D)$ is the rank of the canonical divisor, $g$ is the genus and the $\deg(D)$ is just the degree of the divisor.

## 2.3 Divisors on Smooth Curves and the Projective Line

For this paper, we are only working on curves isomorphic to $\mathbb{P}^1$ (the projective line). The projective line is a convenient way of extending the Euclidean plane, containing points of the form $(x, y)$, to infinity through the use of homogeneous coordinates, which take the form $[x : y]$. This is an important notion because now we are able to speak of the intersection point of every pair of lines, when in the past that was impossible. On the Euclidean plane, parallel lines by definition are lines that never cross. However, on $\mathbb{P}^1$ we simply define the point at infinity as $[1 : 0]$, and say that any two lines that do not cross in the Euclidean plan cross there. A result that immediately follows from this is that $[0 : 0]$ does not exist and that any two sets of coordinates are equivalent up to a scalar. For example, $[1 : 1] = [-3, -3]$, but $[1 : 1] \neq [1 : -1]$. With those facts, we can now define divisors over the smooth curve.

**Definition 2.3.1.** If we are given an irreducible curve X and a function $f\colon X \to \mathbb{P}^1$ we define divisors of curves in this instance as $divf = \sum v_p(f)P$. Here, the coefficients $v_p(f)$ are just the order of the zero or pole of f at P.

These divisors are called principal divisors, which means that they are just divisors of a meromorphic function, or a function that is analytic at all but a certain number of points, with those points going to infinity. It is with this fact that we are able to prove the important result: On a smooth curve, the degree $divf = 0$. This is because for a meremorphic function, there are as many poles as there are zeros. Thus, when we calculate the degree, the zeros take on positive values and the poles take on negative values, leaving us with our desired result. Moving onward, it is important to define the types of maps that will be used by the functions in this paper which are exactly those of the form, $f\colon \mathbb{P}^1 \to \mathbb{P}^1$. These maps take a point $[s : t]$ and send them to a pair of homogeneous polynomials, that is,

$$[s : t] \mapsto [F(s, t) : G(s, t)].$$

Homogeneous polynomials are just those that have every nonzero term of the same degree. Then,

in this instance the **degree** of $f$ is the degree of $F$. But, that means that it is also the number of points in a (generic) fiber, that is, $\deg f = \#f^{-1}([a : b])$. Therefore to solve for the degree we just have to find the values $[s : t]$, that when we plug into $[F(s,t) : G(s,t)]$ the result is the new homogeneous coordinate $[a : b]$. On curves, we define linear equivalence between 2 divisors $D$ and $D'$ if they differ by some principal divisor $div f$. With linear equivalence, we can define the linear system, $|D|$, of $D$ to be

$$\{D' : D' \sim D \text{ and } D' \geq 0\}$$

## 2.4 Divisors on Graph Curves

Before moving on to the main result of this paper, dealing with toric rank, we have one last topic to explain.

**Definition 2.4.1.** A graph curve over $k$, where $k$ is an algebraically closed field, is defined to be a curve $X$ over $k$ such that each component of $X$ is isomorphic to $\mathbb{P}^1$ and such that each pair of components intersects transversely.

Here, "transversely" just means that when two lines meet the equation locally looks like $xy = 0$. Also, we have that the intersection points are called nodes, irreducible curves are those with one component, and a smooth curve is a curve with no nodes.

**Definition 2.4.2.** With a graph curve we are able to associate the dual graph, call it $\Gamma$. The makeup of $\Gamma$ is as follows. The vertices $V = V(\Gamma)$ are the components of $X$ and the edges correspond to the nodes.

We see in Figure 2.2 below, an example of graph curve and its dual graph.

Having that information in hand, we are able to define the idea of a divisor on a graph curve.

**Definition 2.4.3.** A divisor on a graph curve is a formal integral sum $\sum n_P P$ of points, where each $P$ is a smooth point of $X$ and each $nP \in \mathbb{Z}$. With that, we can say $\operatorname{Div} X$ is the free abelian group on the smooth points and $\operatorname{Div} \Gamma$ is the free abelian group on $V(\Gamma)$.

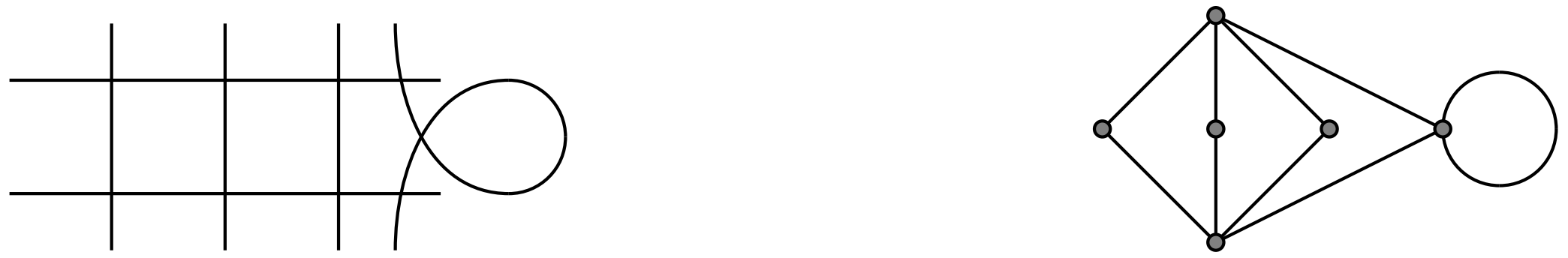

Figure 2.2: A curve and its dual graph.

By labeling the vertices of $\Gamma$, we can now represent $D \in \operatorname{Div}\Gamma$ as a vector. If we do not label the vertices, we refer to a divisor as $\sum n_v v$. The degree of a divisor is just $\sum n_P$. But, this is important because it allows us to create a map from one divisor group to the other. We do this by first labeling the components of $X$ as $X_1, \ldots, X_n$ and labeling the corresponding vertices as $v_1, \ldots, v_n$. Then, the map is just

$$\tau\colon\ \operatorname{Div} X \to \operatorname{Div}\Gamma$$

which sends

$$\sum n_P P \mapsto \sum n_i \left(v_i\right),$$

where

$$n_i = \sum_{P \in X_i} n_P.$$

This also allows us to speak of the specialization of divisors. It is easiest to think of this notion geometrically. The specialization of a divisor is the process by which you collapse the points on a curve to whatever node corresponds to the component they lie on. This is seen in Figure 2.3 below.

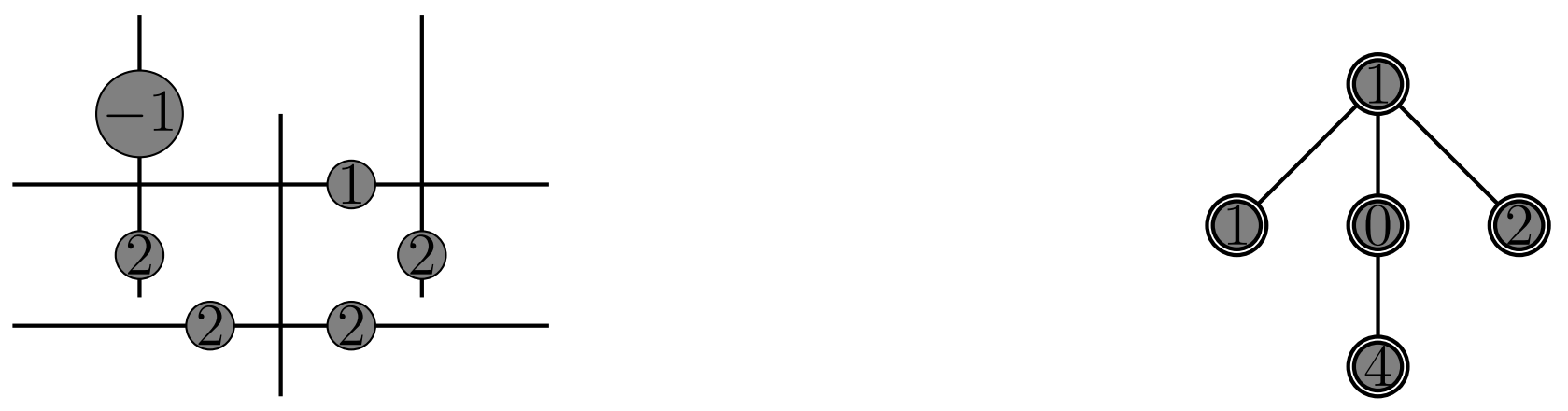

Figure 2.3: Example of specialization on a tree.

# Chapter 3

# Toric Rank

## 3.1 Definition of Toric Rank

Throughout this section, we let $X$ be a graph curve with dual graph $\Gamma$.

**Definition 3.1.1.** Let $X$ be a graph curve. We define the *normalization* $\widetilde{X}$ of $X$ to be the disjoint union $\coprod X_i$ of the irreducible components of $X$. We denote the pre-images of the nodes by $P_{ij} \in X_i$, labeled so that $P_{ij}$ and $P_{ji}$ correspond to the two points above the intersection of $X_i$ and $X_j$. We define a divisor on the normalization to be a divisor on each component curve. Given a divisor $D$ on $X$ we define $\tau(D)$ to be the associated divisor on the dual graph $\Gamma$, which assigns to a vertex associated to a component $X_i$ the sum of the values of $D$ lying on $X_i$.

For simplicity, we will always assume that $P_{ij} \neq \infty$, and this can always be arranged after an automorphism.

**Example 3.1.2.** We have below an example of a curve and its normalization, that have the components labeled, and intersections marked.

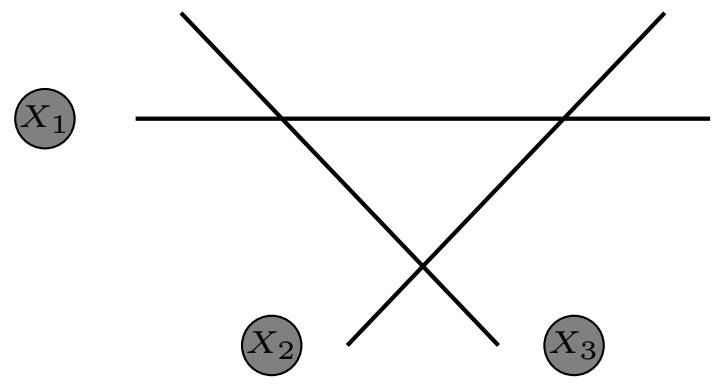

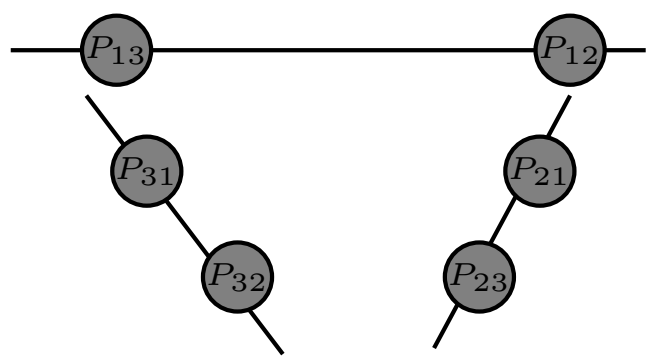

The next important definition to unpack is that of a twisted divisor. A twisted divisor, $D(f)$, is the divisor on the normalization of the graph curve given as follows. Leave each of the points on $D$ as they are. Then, each time we fire, we add and subtract points to the divisor at the intersection points. Formally, this process can be given as follows.

**Definition 3.1.3.** Let $D \in \operatorname{Div} X$ be a divisor and let $f\colon V(\Gamma) \to \mathbb{Z}$ be an integer valued function (representing a sequence of firings). We define $D(f)$ as follows. We have $D(f) = D$ away from the points $P_{ij}$, and given a vertex $v_i$ corresponding to some component $X_i$, we subtract $f(v_i)$ many points from $P_{ij}$ for each $j$ such that $X_i \cap X_j$ is nonempty, and we add $f(v_i)$ points at each such $P_{ji}$.

The next important definition is that of toric rank. First we define what it means to have non-negative rank. Given a non-zero function $f$ on a component $X_i$ of $X$ and for a component $X_j$ that intersects $X_i$, we define $g_j(f)$ to be the function $f(x - P_{ij})^{-v_{P_{ij}}(f)}$. Given a divisor $D$ on the normalization $\widetilde{X}$, we define $D|_{X_i}$ to be the restriction of $D$ to $X_i$.

**Definition 3.1.4.** We say that $r_{\text{tor}}(D) \geq 0$ if there exists a function $f\colon V(\Gamma) \to \mathbb{Z}$ such that

$$D + \operatorname{div} f \geq 0$$

and there exist non-zero functions $f_i \in H^0(X_i, D(f)|X_i)$ such that for every $i, j$ such that $X_i \cap X_j$ is nonempty, we have

$$g_j(f_i)(P_{ij}) = g_i(f_j)(P_{ji}).$$

We say that $r_{\text{tor}}(D) = -1$ if $r_{\text{tor}}(D)$ is not $\geq 0$.

So, colloquially, $r_{\text{tor}}(D) \geq 0$ if $r(\tau(D)) \geq 0$ and, for some equivalence $D \sim D'$ to $D' \in |D|$, witnessed by some firing sequence $f$, the remnants of the twisted divisor $D(f)$ on each component $X_i$ can be made equivalent in a way compatible with the intersection points $P_{ij}$.

Now we come to the general definition, which we define inductively.

**Definition 3.1.5.** We define the toric rank $r_{\text{tor}}(D)$ as follows. We say that $r_{\text{tor}}(D_0) \geq r$ if for every effective divisor $E \in \operatorname{Div} X$ of degree $r$, $r(D_0 - E)_{\text{tor}} \geq 0$.

**Example 3.1.6.** Consider the following tree and dual graph with components (from left to right) $X_1, X_2, X_3$ and intersection points $P_{12}, P_{21}, P_{23}, P_{32}$. For simplicity, assume that the divisor is $1$ at the point $0$ of $X_1$ and $-1$ at the point $0$ of $X_3$. The graph divisor on the right is equivalent to the zero divisor via the firing sequence

$$f(v_1) = 1, f(v_2) = 0, f(v_3) = -1.$$

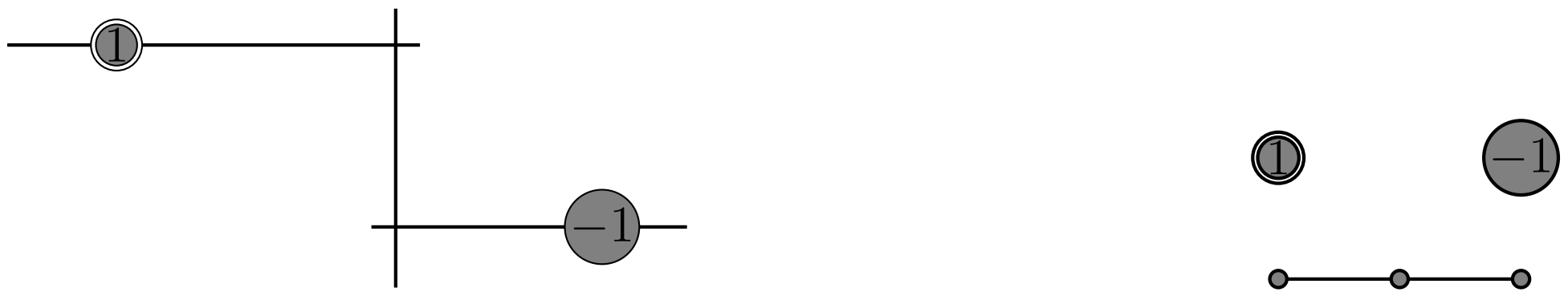

This sequence represents firing at $v_1$ and $v_2$, and results in the zero divisor on the graph. Now, we consider the twisted divisor on the normalization of $X$:

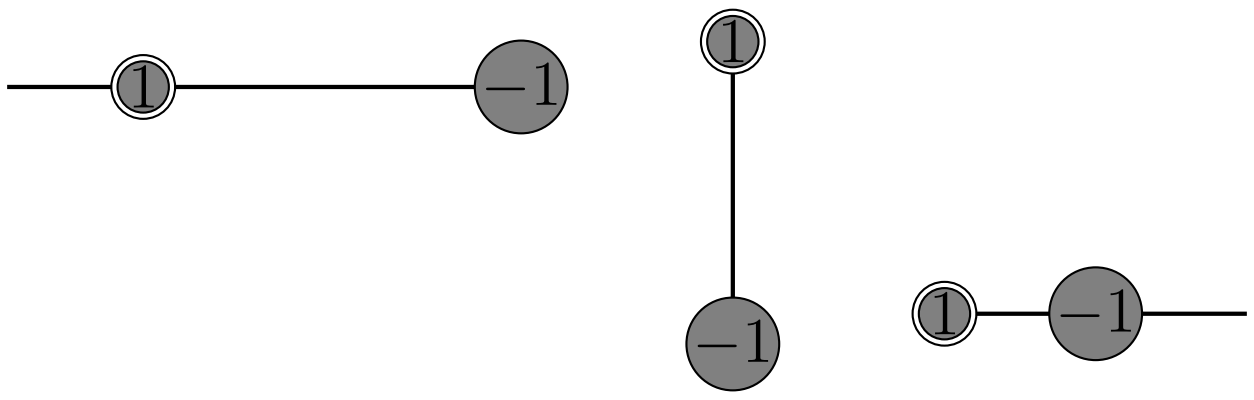

Then the twisted divisor $D(f)$ has degree 0, but moreover the restriction $D(f)|_{X_i}$ has degree zero on each component. It follows that

$$\dim H^0(X_i, D(f)|_{X_i}) = 1$$

for each $i$, so for each $i$ there is a one dimensional space $\langle f_i \rangle = H^0(X_i, D(f)|_{X_i})$ of functions with zeroes and poles constrained by $D(f)$. We can make this explicit as follows, taking

1. $f_1 = A_1(x - P_{12})/x$,

2. $f_2 = A_2(x - P_{23})/(x - P_{21})$,

3. $f_3 = A_3 x/(x - P_{23})$.

To compare these, we modify and get

1. $g_2(f_1) = A_1/x$,

2. $g_1(f_2) = A_2(x - P_{23})$,

3. $g_3(f_2) = A_2/(x - P_{21})$,

4. $g_2(f_3) = A_3 x$.

The compatibility condition is then that

1. $g_2(f_1)(P_{12}) = g_1(f_2)(P_{21})$,

2. $g_3(f_2)(P_{23}) = g_2(f_3)(P_{32})$,

and evaluating gives

$$A_1/P_{12} = A_2(P_{21} - P_{23}) \text{ and } A_2/(P_{32} - P_{21}) = A_3 P_{32}.$$

This is a pair of linear equations, with a solution: $A_2$ is determined by $A_1$, and $A_3$ by $A_2$. We conclude that $r_{\mathrm{tor}}(D) \geq 0$.

## 3.2 Results

**Theorem 3.2.1** (Riemann–Roch Theorem for trees)**.** *Let $X$ be a graph curve whose dual graph is a tree. Then, Riemann–Roch Theorem holds for $X$.*

*Proof.* Let $D \in \operatorname{Div} X$. We may assume that $\deg D \geq 0$, since if $\deg D = -1$ there is nothing to prove, and if $\deg D \leq -2$ we may replace $D$ with $D - K$, which has non-negative degree. Since $\deg D \geq 0$ implies that $\deg K - D \leq 0$, $r_{\mathrm{tor}}(K - D) = -1$. Thus, the goal is to prove that

$$r_{\mathrm{tor}}(D) = \deg D.$$

For this, it suffices (from the definition of rank) to prove that if $\deg D$ is exactly zero, then it is torically equivalent to an effective divisor.

Now, Riemann–Roch Theorem holds for $\tau(D)$, so there exists a sequence $f\colon V(\Gamma) \to \mathbb{Z}$ of chip firings such that $D(f)|_{X_i}$ has non-negative degree. On each component, the space $H^0(X_i, D(f)|_{X_i})$ has dimension at least one. We proceed by a "depth-first search". Let $X_1$ be any component, and let $f_1$ be any non-zero element of $H^0(X_1, D(f)|_{X_1})$. Consider the components intersecting $X_1$ non-trivially; since $X$ is a tree, they intersect exactly once; the intersection with $X_1$ imposes exactly one condition on each $f_j$. For such a component $X_j$, since $\dim H^0(X_j, D(f)|_{X_j}) \geq 1$, we can choose $f_j$ so that $g_i(f_1)(P_{1j}) = g_1(f_j)(P_{ji})$. Next, a component $X_k$ intersecting $X_1$ or the previous $X_j$, intersects at most once, so we can again arrange compatibility of the elements $f_k$. Since $\Gamma$ is a tree, as we proceed, no component ever intersects the part of the tree we have dealt with, and so by induction we can choose non-zero $f_i$ compatibly on each component, completing the proof. □

**Lemma 3.2.2.** *Let $\Gamma' \subset \Gamma$ be a subgraph and suppose that the complement is a union of trees such that each tree intersects $\Gamma$ at most once. Then, Riemann–Roch Theorem holds for $\Gamma$ if and only if it holds for $\Gamma'$.*

*Proof.* The proof is the same as that of Theorem 3.2.1. □

**Theorem 3.2.3** (Riemann–Roch Theorem for $g(\Gamma) = 1$)**.** *Let $X$ be a graph curve whose dual graph has genus 1. Then, Riemann–Roch Theorem holds for $X$.*

*Proof.* Let $D \in \operatorname{Div} X$. First, since $\deg K = 0$, if $\deg D \leq -1$ then by symmetry it suffices to prove Riemann–Roch Theorem for $K - D$. Moreover, if $\deg D = 0$, then since $K = 0$, $K - D = -D$, and in this case $r_{\text{tor}}(D) = r_{tor}(-D)$, so Riemann–Roch Theorem holds by default (since both sides are 0).

We may assume that $\deg D \geq 0$, in which case the goal is to prove that

$$r_{\text{tor}}(D) = \deg D - 1.$$

For this, it suffices (from the definition of rank) to prove that if $\deg D$ is exactly one, then it is torically equivalent to an effective divisor.

Now, by Lemma 3.2.2 we may assume that $\Gamma$ is an $n$-gon. Riemann–Roch Theorem holds for $\tau(D)$, so there exists a sequence $f\colon V(\Gamma) \to \mathbb{Z}$ of chip firings such that $D(f)|_{X_i}$ has non-negative degree, and at least one component (lets renumber so that it is the first component) $D(f)|_{X_1}$ has degree one. The spaces $H^0(X_i, D(f)|_{X_i})$ have dimension one for $i \neq 1$ and two for $i = 1$. The complement of $X_1$ is a tree, so as in the proof of Theorem 3.2.1 we may arrange the existence of non-zero functions $f_i \in H^0(X_i, D(f)|_{X_i})$ (for $i \neq 1$) with agreement at the nodes. Since $\dim H^0(X_1, D(f)|_{X_1}) = 2$, we may arrange for a function $f_1$ with agrees at the nodes too, completing the proof.

□

Next, we note that Riemann–Roch Theorem predicts that

$$r_{\text{tor}}(D) = r_{\text{tor}}(K - D) + \deg D + 1 - g \geq \deg D - g$$

(since $r_{\text{tor}}(K - D) \geq -1$). In particular, if Riemann–Roch is true, then if $\deg D \geq g$, then $r_{\text{tor}}(D) \geq 0$ and $D$ is torically equivalent to an effective divisor. We prove an easier corollary of Riemann–Roch, which we refer to as "asymptotic Riemann–Roch".

**Theorem 3.2.4** (Asymptotic Riemann–Roch Theorem)**.** *If* $\deg D \geq 2g$*, then* $r_{\text{tor}}(D) \geq 0$.

*Proof.* We proceed by another "depth search". We define a sequence of subgraphs $\Gamma_1, \ldots, \Gamma_n$ as follows. We let $v_1$ be any vertex and define $\Gamma_1$ to be just $v_1$. For $\Gamma_i$, we pick any vertex $v_i$ such that $v_i \notin \Gamma_{i-1}$, but such that there is an edge from $v_i$ to at least one vertex of $\Gamma_{i-1}$, and we define $\Gamma_i$ to be the maximal subgraph containing $\Gamma_{i-1}$ and $v_i$.

Define the divisor $D'$ to be $D' := \sum_i \max\{\deg_{\Gamma_i}(v_i), 0\}(v_i)$, where $\deg_{\Gamma_i}(v_i)$ is the degree of $v_i$ as a vertex of $\Gamma_i$. (Equivalently, $\deg_{\Gamma_i}(v_i)$ is the number of edges from $v_i$ to $\Gamma_{i-1}$.) We claim that $\deg D' = g$. Indeed, it is clear from the Euler characteristic formula that each edge from $v_i$ to $\Gamma_{i-1}$ beyond the first increases the genus by 1, so that $g(\Gamma_i) = g(\Gamma_{i-1}) + \max\{\deg_{\Gamma_i}(v_i), 0\}$.

Now, by Riemann–Roch for graphs,

$$r(D-D') \geq \deg(D-D') - g \geq 2g-2-2 \geq 0,$$

so there exists a sequence $f\colon V \to \mathbb{Z}$ such that $D - D' + \Delta(f) \geq 0$, i.e. $D + \Delta(f) = D' + E$ for some effective divisor $E$. To show that $r_{\mathrm{tor}}(D) \geq 0$, we need to show that the twisted divisor $D(f)$ is torically equivalent to an effective divisor. By our definition of $D'$ and effectivity of $E$, we have

$$\deg D(f)|_{X_i} \geq \deg_{\Gamma_i}(v_i),$$

and in particular

$$\dim H^0(X_i, D(f)|_{X_i}) \geq \deg_{\Gamma_i}(v_i)$$

We proceed inductively. We let $f_1 \in H^0(X_1, D(f)|_{X_1})$ be any non-zero function. The intersection of $v_2$ with $\Gamma_1$ imposes at most $\deg_{\Gamma_2}(v_2)$ conditions that a function in $f_1 \in H^0(X_2, D(f)|_{X_2})$ must satisfy, and we can find such a function since $\dim H^0(X_2, D(f)|_{X_2}) \geq \deg_{\Gamma_2}(v_2)$. Similarly, the intersection of $v_i$ with $\Gamma_{i-1}$ imposes at most $\deg_{\Gamma_i}(v_1)$, and since $\dim H^0(X_i, D(f)|_{X_i}) \geq \deg_{\Gamma_i}(v_i)$, there exist a compatible function $f_i \in H^0(X_i, D(f)|_{X_i})$. By induction, we can choose compatible $f_i$, and we are done.

□

To proceed, using a more delicate linear algebra argument we prove a stronger asymptotic theorem.

**Theorem 3.2.5** (Stronger Asymptotic Riemann–Roch Theorem)**.** *If* $\deg D \geq g$*, then* $r_{\mathrm{tor}}(D) \geq 0$.

*Proof.* Label the vertices $v_1, \ldots, v_n$. Let $D(f)$ be a twisted divisor such that $\deg D(f)_{X_i} \geq 0$ for all $i$ and set $k_i := H^0(X_i, D(f)|_{X_i}) = \deg D(f)_{X_i} + 1$. Given a function $f_i \in H^0(X_i, D(f)|_{X_i})$, define the column vector $T(f_i)$ to have $j$th component 0 if $v_i$ is not adjacent to $v_i$ and $g_j(f_i)(P_{ij})$ otherwise. For each $i$, let $f_{i,1}, \ldots, f_{i,k_i}$ be a basis for $f_i \in H^0(X_i, D(f)|_{X_i})$. Then the vectors $T(f_{i,1}), \ldots, T(f_{i,k_i})$ are linearly independent.

Let $M$ be the matrix whose columns are $T(f_{i,j})$ (where both $i$ and $j$ vary). Then the agreement condition can be phrased as the existence of a vector $\overline{v}$ such that $M\overline{v} = 0$ (since the rows summing to zero implies that there is a linear combination of functions at each component which agree) and such that the entries of $\overline{v}$ corresponding to the coefficients of the functions $f_{i,j}$ for *fixed* $i$ are not all zero (this, together with the linear independence in the last paragraph, ensures that the resulting functions are non-zero on each component).

Now a counting argument will ensure the existence of a non-trivial null space as follows. The number of rows is equal to the number of edges, and the number of columns is $\deg D(f) + g$. Thus, if $\deg D(f) + g$ is greater than the number of edges, the corresponding matrix will have more columns than rows and thus a non-trivial vector in the kernel. Linear independence of the $f_{i,k}$'s will ensure that we can choose this vector to have non-zero entries, completing the proof.

□

**Example 3.2.6.** Consider the following graph curve

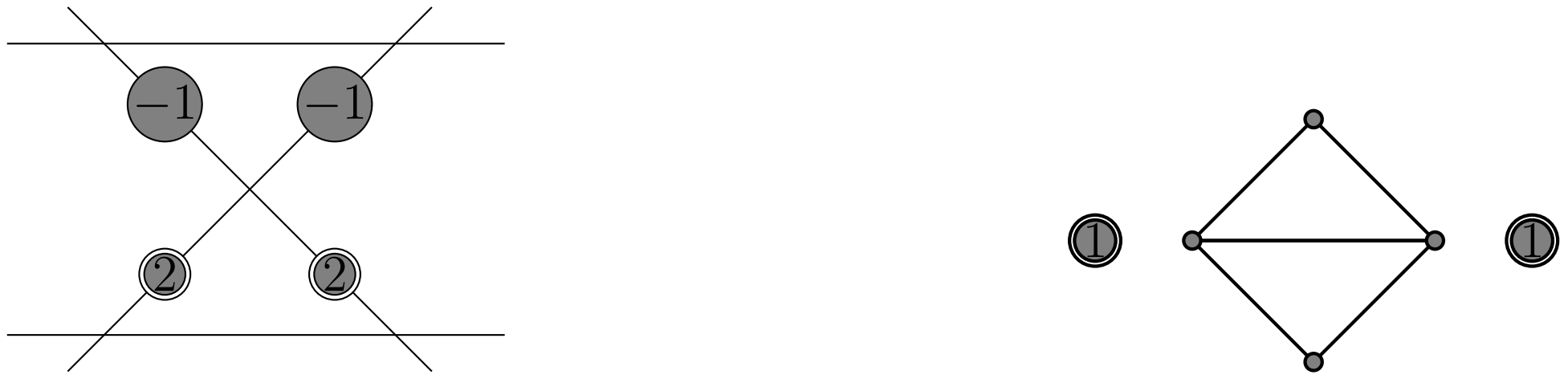

The associated $D$ divisor on the graph is already effective, so there is no twisting to do. Riemann–Roch predicts that $D$ is torically equivalent to an effective divisor. Label the components of $X$ clockwise, beginning with the top horizontal curve. Then, for instance, on $X_i$, $D|_{X_1}$ has degree 0, so $H^0(X_1, D|_{X_1})$ has dimension 1. The rows correspond to the conditions at the nodes $P_{12}$, $P_{13}$, $P_{13}$, $P_{13}$, $P_{13}$, giving us the matrix

$$\begin{pmatrix} * & * & * & 0 & 0 & 0 \\ * & 0 & 0 & * & * & 0 \\ 0 & * & * & * & * & 0 \\ 0 & * & * & 0 & 0 & * \\ 0 & 0 & 0 & * & * & * \end{pmatrix}$$

Since the rank of this matrix is 5, there is a non-trivial vector in the null space, which specifies the correct linear combination of functions to take.

# Chapter 4

# Experimental Data

## 4.1 Data

Using the Magma code in the next section, we tested Riemann–Roch for a large space of divisors. For low genus, we enumerated all "tree-less" graphs of small genus and tested Riemann–Roch for all divisors of degree at most $g - 1$. For larger genus, we picked random graphs and tested Riemann–Roch for all divisors on such graphs. This becomes difficult to do exhaustively around genus $4$ and impossible to do exhaustively around genus 6. For larger genus, we pick a random divisor and test Riemann–Roch for this single divisor. In all cases, Riemann–Roch is found to hold.

## 4.2 Code

```
// Checks whether the graph with given adjacency matrix is connected
isConnected := function(M)
  if not IsSymmetric(M) then
    return false;
  else
    return &and{N[i][j] ge 1 : i,j in [1..Nrows(M)] }
```

```
        where N is &+[M^i : i in [1..Nrows(M)]];
  end if;
end function;

randomAdjacencyMatrix := function(n)
  M:= Matrix(n,n,[[ i ge j select 0 else
            Random([0,1]) : i in [1..n]] : j in [1..n]]);
  return M + Transpose(M);
end function;

randomConnectedAdjacencyMatrix := function(n)
  M := randomAdjacencyMatrix(n);
  while not isConnected(M) do
  M := randomAdjacencyMatrix(n);
  end while;
  return M;
end function;

canonicalDivisor := function(M)
  return Matrix(Nrows(M),1,[ &+[M[i][j]
        : j in [1..Nrows(M)]] - 2 : i in [1..Nrows(M)]]);
end function;

degree := function(D)
  return (&+[ D[i] : i in [1..Nrows(D)]])[1]; end function;
genus := function(M);
  return Integers()!((degree(canonicalDivisor(M))+2)/2); end function;
```

```
produceNonEffectiveDivisorsOfDegree := function(v,n,m)
  // v := (integer) number of verticies
  // n := (integer) sum of chips
    return {@
     [a : a in tup] : tup in CartesianPower([-m..n+m],v)
     | &+[a : a in tup] eq n @};
end function;

produceDivisorsOfDegree := function(v,n)
  // v := (integer) number of verticies
  // n := (integer) sum of chips
    return {@
     [a : a in tup] : tup in CartesianPower([0..n],v)
   | &+[a : a in tup] eq n
    @};
end function;

adjacencyToLaplacian := function(M)
n := Degree(Parent(M));
A := M*Matrix(n,1,[1 : i in [1..n]]);
D := DiagonalMatrix([ A[i][1] : i in [1..n]]);
return -1*M + D;
end function;

toricTest := function(M,c : B := 10^10)
  N := (Matrix(
```

```
  [
    [not l in {i,j} select 0 else Random(FiniteField(NextPrime(B)))
     : m in [1..c[l][1]+1],  l in [1..Nrows(c)*Ncols(c)]
    ]
     : k in [1..M[i][j]]  , j in [i+1..Nrows(M)], i in [1..Nrows(M)]
  ]
  ));
  v := Random(Nullspace(Transpose(N)));
  return not 0 in {v[i] : i in [1..Ncols(v)]}, N;
end function;



BoundingBoxToMatrix := function(R,n)
     dual := Dual(Ambient(R));
     return Matrix([ [dual.i*Vertices(R)[j]
            : i in [1..n]] : j in [1..#Vertices(R)] ] );
end function;


  findMin:= function(M)
  M:= Transpose(M);
  min:= [];
  for i in [1..Nrows(M)] do
     min[i]:= Minimum([M[i][j] : j in [1..Ncols(M)]]);
  end for;
  return min;
end function;
```

```
findMax:= function(M)
M:= Transpose(M);
max:= [];
for i in [1..Nrows(M)] do
     max[i]:= Maximum([M[i][j] : j in [1..Ncols(M)]]);
end for;
return max;
end function;

getTuples:= function(min,max)
tups:= [];
for i in [1..#min] do;
   tups[i]:= {a : a in [min[i]..max[i]]};
end for;
return tups;
end function;

produceOnlyIntegerSolutions :=function(Q)
for i in [1..Nrows(Q)] do
    for j in [1..Ncols(Q)] do
   if Q[i][j] ge 0 then
      Q[i][j] := (Floor(Q[i][j]));
   else
      Q[i][j] := (Ceiling(Q[i][j]));
  end if;
    end for;
end for;
```

```
return Q;
end function;

if degree(c) le -1 then return {}; else

v:= Rank(Parent(c));   // number of verticies
A:= adjacencyToLaplacian(M);
D := { }; // this will contain the linear system
R := [] ;

for i in [1 .. v] do
 Mtemp := RemoveColumn(A,  i);
P := PolyhedronWithInequalities([[Mtemp[k][j] : j in [1..v-1]]
     : k in [1..v]] ,[-c[k][1] : k in [1..v]]);
R[i] := BoundingBox(P);
end for;



Q:= [produceOnlyIntegerSolutions(
     BoundingBoxToMatrix(R[i],v-1)) : i in [1..#R] ];

min :=  [findMin(Q[i]) : i in [1..#Q]];
max := [findMax(Q[i]) : i in [1..#Q]];
tups := [getTuples(min[i],max[i])  : i in [1..#min]];
tups := [Insert(tups[i],i,{0}) : i in [1..#tups]  ];
```

```
C:= [CartesianProduct(tups[i]) : i in [1..#tups]];

for i in [1..#C] do
for tup in C[i] do
  tempTup := Matrix(#tup,1,[a : a in tup]);
   if  &and{(c+A*tempTup)[i][1] ge 0 : i in [1..v]} then
     D:= D join {c+A*tempTup};
   end if;
end for;
end for;

return D;
end if;
end function;

produceRank  := function(M,c)
  v := Nrows(c)*Ncols(c);
  linSys := linearSystem(M,c);
  r := -1;
  rankKnown := false;

  while not rankKnown do
  E := produceDivisorsOfDegree(v,r+1);
  rankKnown := exists{e : e in E |
   not exists{
     d : d in linSys |
     &and{(d - Matrix(v,1,e))[i][1]  ge 0 : i in [1..v]} }
```

```
                         };
  if not rankKnown then
    r := r+1;
  end if;
  end while;
  return r;
end function;

produceToricRank  := function(M,c)
  v := Nrows(c)*Ncols(c);
  linSys := linearSystem(M,c);
  r := -1;
  rankKnown := false;

  while not rankKnown do
  E := produceDivisorsOfDegree(v,r+1);
  rankKnown := exists{e : e in E |
    not exists{d : d in linSys |
       &and{(d - Matrix(v,1,e))[i][1]  ge 0 : i in [1..v]}
                         and toricTest(M,d-Matrix(v,1,e))
                         }};
  if not rankKnown then
    r := r+1;
  end if;
  end while;
  return r;
end function;
```

```
singleTestRR := function(M,c)
  K := canonicalDivisor(M);
  return produceRank(M,c) - produceRank(M,K-c) - degree(c)  - 1
         + (degree(K) + 2)/2 eq 0;
  // note that deg K = 2g-2
end function;

singleToricTestRR := function(M,c)
  K := canonicalDivisor(M);
  return produceToricRank(M,c) - produceToricRank(M,K-c)
         - degree(c)  - 1 + genus(M) eq 0;
  // note that deg K = 2g-2
end function;

/////////////////////////////////////////////////////////////
// testing
// Now, a big loop to test RR
v := 7;
while true do
  i := -1;
    M := randomConnectedAdjacencyMatrix(v);
    g := genus(M);
//if g eq 1 then
"**********************************************************";
    "genus",  g;//M;
```

```
    while i le g-2 do
//    while i le 2*g do // we know RR holds for i geq g.
                        // Switch these lines to test this.
      i := i + 1;
      "degree", i;
//      for c in produceNonEffectiveDivisorsOfDegree(v,0,i) do
        for c in produceDivisorsOfDegree(v,i) do
        if not singleToricTestRR(M,Matrix(v,1,c))
           then a := c;MM := M;c; end if;
      end for;
  end while;
//end if;
end while;

/////////////////////////////////////////////////////////////
// testing
// Test RR for a few random graphs of large genus
while true do
  v := Random([5..10]);
    M := randomConnectedAdjacencyMatrix(v);
    g := genus(M);
if g ge 4 then
"**********************************************************";
    "genus",  g;//M;
      i := g-1;
      "degree", i;
        c := Random(produceDivisorsOfDegree(v,i));
```

```
        if not singleToricTestRR(M,Matrix(v,1,c))
           then a := c;MM := M;c; end if;
end if;
end while;
```